\documentclass{amsart}

\usepackage{amssymb,amsthm,mathrsfs}
\usepackage{hyperref}

\newcommand{\arxiv}[1]{arXiv:\linebreak[0]\href{http://arxiv.org/abs/#1}{#1}}

\newtheorem*{thm}{Theorem}
\theoremstyle{definition}
\newtheorem*{prob}{Problem}

\begin{document}

\title[The independence number of the orthogonality graph]{The Independence Number of the Orthogonality Graph in Dimension $2^k$}

\author{Ferdinand Ihringer}
\address{Department of Mathematics: Analysis, Logic and Discrete Mathematics, Ghent University, Belgium}
\thanks{The first author is supported by a postdoctoral fellowship of the Research Foundation --- Flanders (FWO)}
\email{ferdinand.ihringer@ugent.be}

\author{Hajime Tanaka}
\address{\href{http://www.math.is.tohoku.ac.jp/}{Research Center for Pure and Applied Mathematics}, Graduate School of Information Sciences, Tohoku University, Japan}
\thanks{The second author is supported by JSPS KAKENHI Grant Number JP17K05156.}
\email{htanaka@tohoku.ac.jp}

\begin{abstract}
  We determine the independence number of the orthogonality graph on $2^k$-dimensional hypercubes.
  This answers a question by Galliard from 2001 which is motivated by a problem in quantum information theory.
  Our method is a modification of a rank argument due to Frankl 
  who showed the analogous result for $4p^k$-dimensional hypercubes, where $p$ is an odd prime.
\end{abstract}

\maketitle

\section{Introduction}

The \textit{orthogonality graph} $\Omega_n$ has the elements of $\{ -1, 1 \}^n$ as vertices, and two vertices
are adjacent if they are orthogonal, in other words, if their Hamming distance is $n/2$.
The graph $\Omega_n$ occurs naturally when comparing classical and quantum communication \cite{BCW1998}.
In particular, for $n=2^k$ the cost of simulating a specific quantum entanglement on $k$ qubits 
can be reduced to determining the chromatic number $\chi(\Omega_{n})$ of $\Omega_n$ \cite{BCT1999,Galliard2001D}.
The graph $\Omega_n$ is edgeless if $n$ is odd, and is bipartite if $n\equiv 2 \pmod{4}$.
For $n\equiv 0 \pmod{4}$, Frankl \cite{Frankl1986C} and Galliard \cite{Galliard2001D} constructed an independent set of $\Omega_n$ of size
\begin{equation*}\label{pin}
	a_n := 4 \sum_{i=0}^{n/4-1} \binom{n-1}{i},
\end{equation*}
and Galliard \cite{Galliard2001D} asked in 2001 if this is the independence number $\alpha(\Omega_n)$ of $\Omega_n$ when $n=2^k$, $k\geqslant 2$.
Newman \cite{Newman2004D} and, according to \cite[p.~275, Remark]{FR1987}, Frankl conjectured that this holds whenever $n\equiv 0\pmod{4}$.
See also \cite{Cameron2010EJC}. 
Frankl \cite{Frankl1986C} already showed the conjecture in 1986 for all $n=4p^k$ for $k \geqslant 1$, where $p$ is an odd prime.
De Klerk and Pasechnik \cite{KP2007EJC} proved the conjecture for $n=16$, i.e., that $\alpha(\Omega_{16})=2304$, 
using Schrijver's semidefinite programming bound \cite{Schrijver2005IEEE}. Furthermore,
Frankl and R\"{o}dl \cite{FR1987} showed that $\alpha(\Omega_{n}) < 1.99^n$ if $n\equiv 0\pmod{4}$.
In this note, we apply Frankl's method from \cite{Frankl1986C} to show the following:
\begin{thm}
Let $n = 2^k$ for some $k \geqslant 2$.
Then $\alpha(\Omega_n)=a_n$.
\end{thm}
Together with the discussion in \cite[Section 5.5]{Galliard2001D}, that is using $\chi(\Omega_n) \geqslant 2^n/\alpha(\Omega_n)$, 
our result implies an explicit version of Theorem 4 in \cite{BCT1999}.
Finding such an explicit result is one motivation for Galliard's work.
See also \cite{GTW2003P,GN2008SIAM}.

\section{Proof of the Theorem}

Let $A_j$ be the $0$-$1$-matrix indexed by the vertices of the hypercube $Q_n = \{ -1, 1 \}^n$ with $(A_j)_{xy} = 1$ if 
$x$ and $y$ have Hamming distance $j$.
The matrices $A_j$ have $n+1$ common  eigenspaces $V_0, V_1, \dots, V_n$, and
in the usual ordering of the eigenspaces the eigenvalue of $A_j$ with respect to $V_i$
is given by the Krawtchouk polynomial (see \cite[Theorem 4.2]{Delsarte1973})
\begin{align*}
 	K_j(i)=K_j(i; n) := \sum_{h=0}^j (-1)^h \binom{i}{h} \binom{n-i}{j-h}.
\end{align*}
It is known that the orthogonal projection matrix $E_i$ onto $V_i$ has
the entry $(E_i)_{xy} = 2^{-n} K_i(j)$ if $x$ and $y$ are at Hamming distance $j$ \cite[Theorem 4.2]{Delsarte1973}, so that we have in particular $\operatorname{rank} E_i=\operatorname{trace} E_i=K_i(0)= \binom{n}{i}$.
The $(n+1)$-dimensional matrix algebra
spanned by $A_0=I, A_1, \dots, A_n$ is called the \textit{Bose--Mesner algebra} of $Q_n$.

Assume now that $n=2^k$, $k\geqslant 3$.
(The result is trivial if $k=2$.)
Let $C$ be an independent set of $\Omega_{2^k}$, and let $C_{\mathrm{even}}^{\pm},C_{\mathrm{odd}}^{\pm} \subseteq \{ -1, 1 \}^{2^k-1}$ be as in \cite{Frankl1986C}:
$C_{\mathrm{even}}^+$ is given by taking all the even-weight elements of $C$ that end with $+1$, followed by truncating at the last coordinate, and the other three are analogous.
Let $C'$ be one of these four families.
Then the Hamming distances in $C'$ are even and unequal to $2^{k-1}$, so they lie in the following set:
\begin{equation}\label{weights}
	\{2s : s=0,1,\dots,2^{k-1}-1,\, s\ne 2^{k-2}\}.
\end{equation}

Below we work with the Bose--Mesner algebra $\mathscr{A}$ of $Q_{2^k-1}$.
For every $M\in\mathscr{A}$, let $\overline{M}$ denote the principal submatrix corresponding to $C'$.
Consider the polynomial
\begin{equation*}
	\varphi(\xi)=\binom{\xi/2-1}{2^{k-2}-1} \in \mathbb{R}[\xi],
\end{equation*}
and expand it in terms of the Krawtchouk polynomials $K_i(\xi)=K_i(\xi;2^k-1)$:
\begin{equation}\label{Krawtchouk expansion}
	\varphi(\xi)=\sum_{i=0}^{2^{k-2}-1} c_i K_i(\xi).
\end{equation}
Let
\begin{equation*}
	X=\sum_{j=0}^{2^k-1}\varphi(j)A_j \in\mathscr{A}.
\end{equation*}
On the one hand, observe that $\overline{X}$ has only integral entries in view of \eqref{weights}, and an easy application of Lucas' theorem (cf.~\cite{Fine1947AMM}) shows moreover that $\overline{X}\equiv \overline{I} \pmod{2}$.
In particular, $\overline{X}$ is invertible.
On the other hand, from \eqref{Krawtchouk expansion} we have
\begin{equation*}
	X=2^{2^k-1}\sum_{i=0}^{2^{k-2}-1} c_i E_i.
\end{equation*}
It follows that
\begin{align*}
	|C'| &= \operatorname{rank} \overline{X} \leqslant \operatorname{rank} X 
	 \leqslant \sum_{i=0}^{2^{k-2}-1} \operatorname{rank} E_i = \sum_{i=0}^{2^{k-2}-1} \binom{2^k-1}{i}.
\end{align*}
As $|C| = |C_{\mathrm{even}}^+| + |C_{\mathrm{even}}^-| + |C_{\mathrm{odd}}^+| + |C_{\mathrm{odd}}^-|$, the theorem follows.

\section{Future Work}

Schrijver's semidefinite programming bound has been extended to hierarchies of upper bounds; see, e.g., \cite{BGSV2012B,Laurent2007MP}.
In view of \cite{KP2007EJC}, it is interesting to investigate if these bounds in turn prove the conjecture for other values of $n$.
One of the referees pointed out to us that using next level in the hierarchy, see \cite{GMS2012IEEE}, yields the correct bound
of $a_{24}=178208$ for the case $n=24$.

\begin{prob}
Prove the conjecture for $n=40$, which is the first open case.
\end{prob}

\paragraph*{Acknowledgements} We thank the anonymous referee for solving the case $n=24$.

\end{document}